\documentclass[11pt]{article}
\usepackage{latexsym}
\usepackage{amssymb}

\def\qed{\hfill$\Box$}

\newcommand{\bydef}{\mbox{$ \;\stackrel{\triangle}{=}\; $}}

\newcommand{\subD}{_{{}_D}}

\newcommand{\RL}{{\mathbb R}}

\newcommand{\IND}{{\mathbb I}}

\newcommand{\PR}{\mbox{\rm Pr}}

\newcommand{\Ahat}{\hat{A}}

\newcommand{\calE}{{\cal E}}

\newcommand{\Dmax}{\mbox{$D_{\rm max}$}}

\newcommand{\Rmax}{\mbox{$R_{\rm max}$}}

\def\be{\begin{eqnarray}}
\def\ee{\end{eqnarray}}
\def\ben{\begin{eqnarray*}}
\def\een{\end{eqnarray*}}

\input epsf

\title{A Remark on Unified Error Exponents:\\
Hypothesis Testing,
Data Compression \& Measure Concentration}

\author{I. Kontoyiannis \and A.D. Sezer}
 
\date{\today}
 
\begin{document}
\bibliographystyle{plain}
\maketitle
 
\thispagestyle{empty}

\footnotetext[1]{
Ioannis Kontoyiannis is with the Division
of Applied Mathematics and the Department
of Computer Science, Brown University,
Box F, 182 George St., Providence, RI 02912, USA.
Email:~{\tt yiannis@dam.brown.edu}
Web:~{\tt www.dam.brown.edu/people/yiannis/}
}
 
\footnotetext[2]{
Ali Devin Sezer 
is with the Division
of Applied Mathematics, Brown University,
Box F, 182 George St., Providence, RI 02912, USA.
Email:~{\tt ali\_sezer@brown.edu}
}
 
\footnotetext[3]{
I.\ Kontoyiannis was supported in part
by NSF grants \#0073378-CCR
and DMS-9615444, and by 
USDA-IFAFS grant \#00-52100-9615.}

\vspace{0.3in}

\centerline{\bf Abstract}
{\small
\begin{quote}
Let $A$ be finite set equipped
with a probability distribution $P$,
and let $M$ be a ``mass'' function on $A$.
A characterization is given for the most 
efficient way in which $A^n$ can be 
covered using spheres of a fixed radius. 
A covering is a subset $C_n$ of $A^n$
with the property that most of the elements 
of $A^n$ are within some fixed distance 
from at least one element of $C_n$, and
``most of the elements''
means a set whose probability 
is exponentially close to one
(with respect to the product 
distribution $P^n$). An efficient covering 
is one with small mass~$M^n(C_n)$.
With different 
choices for the geometry on $A$,
this characterization gives various 
corollaries as special cases,
including Marton's error-exponents theorem in
lossy data compression,
Hoeffding's optimal hypothesis testing 
exponents, and a new sharp converse 
to some measure concentration 
inequalities on discrete spaces. 
\end{quote}

\vspace{0.4in}


}

\newpage 

\section{Introduction}

Let $A$ be a finite set and
$P$ a probability distribution on $A$.
Suppose that the distance (or ``distortion'')
$\rho(x,y)$ between any two points 
$x,y\in A$ is measured by a given 
nonnegative function
$\rho:A\!\times\!A\to[0,\infty)$,
and for strings
$x_1^n=(x_1,x_2,\ldots,x_n)$ and
$y_1^n=(y_1,y_2,\ldots,y_n)$ in $A^n$
let $\rho_n(x_1^n,y_1^n)$ be the
corresponding coordinate-wise distance
(or single-letter distortion measure) on $A^n\times A^n$:
\ben
\rho_n(x_1^n,y_1^n)=\frac{1}{n}\sum_{i=1}^n\rho(x_i,y_i).
\label{eq:sldistance}
\een
Since $A$ is a finite set, the function $\rho$
is bounded above by 
$$\Dmax\bydef \max_{x,y\in A}\rho(x,y)
=\max_{x_1^n,y_1^n\in A^n}\rho_n(x_1^n,y_1^n).$$

Without loss of generality
we assume throughout 
that $P(a)>0\,$ for all $a\in A$, and
that for each $\,a\in A$ there exists
a $\,b\in A$ with $\rho(a,b)=0$
(otherwise
we may consider 
$\rho'(x,y)=[\rho(x,y)-\min_{z\in A}\rho(x,z)]$
instead of $\rho(x,y)$).

Given a $D\geq 0$, we want to cover 
``most'' of $A^n$ using balls $B(y_1^n,D)$,
where
\ben
B(y_1^n,D)=\{x_1^n\in A^n\;:\;\rho_n(x_1^n,y_1^n)\leq D\}
\een
is the closed ball of radius $D$ centered 
at $y_1^n\in A^n$. 
To be precise, given a set $C_n\subset A^n$,
we write $[C_n]\subD$ for the 
{\em $D$-blowup of $C_n$},
$$[C_n]\subD\bydef \bigcup_{y_1^n\in C_n}
        \!\!\!B(y_1^n,D).$$
A {\em $D$-covering of $A^n$}
is a sequence of subsets $C_n$ of
$A^n$, $n\geq 1$, such that 
the $P^n$-probability of the part of 
$A^n$ which is not covered by $C_n$ 
within distance $D$
has exponentially small probability,
\be
\PR\{\mbox{``error''}\}
\bydef
	1-P^n([C_n]\subD)
&\approx&
	2^{-nE},
\label{eq:cover}
\ee
for some $E>0$.
We are interested in ``efficient'' coverings
of $A^n$, that is, given a ``mass function''
$M:A\to(0,\infty)$, we want to find $D$-coverings
$\{C_n\}$ that 
satisfy (\ref{eq:cover})
and also have small mass
\ben
M^n(C_n)\bydef \sum_{y_1^n\in C_n} M^n(y_1^n)
=\sum_{y_1^n\in C_n}\prod_{i=1}^nM(y_i).
\een
Clearly there is a trade-off between 
finding coverings $\{C_n\}$ with small mass, 
and coverings with a good (i.e., large) 
error-exponent $E$ as in (\ref{eq:cover}). 
Typically, the better the error-exponent, 
the larger the $C_n$,
and the bigger their mass would tend to be.

Motivated, in part, by the following example
and by the applications
illustrated in the examples of the
following section, in our main result 
we give a precise characterization of 
this trade-off.

\medskip
 
\noindent
{\em Example: Measure Concentration on the Binary Cube. }
 
Consider the $n$-dimensional binary cube
$A^n=\{0,1\}^n.$ We measure distance
on $A^n$ by the
proportion of mismatches between two binary
strings $x_1^n$ and $y_1^n$, i.e., we take
$\rho_n(x_1^n,y_1^n)$ to be the Hamming distance,
\be
\rho_n(x_1^n,y_1^n)=\frac{1}{n}\sum_{i=1}^n
\IND_{
\{x_i\neq y_i\}
},\;\;\;\;
x_1^n, y_1^n\in A^n,
\label{eq:hammingdist}
\ee
which also coincides with
the normalized graph distance
when $A^n$ is
equipped with the nearest-neighbor
graph structure. For simplicity,
in this example we consider
natural logarithms
and exponentials.
 
A well-known measure concentration
inequality
\cite[Prop.~2.1.1]{talagrand:95}\cite[Thm.~3.5]{mcdiarmid:98}
gives a precise lower bound on the
sphere-covering error probability
of an arbitrary $C_n$:
For any $D\geq 0$, any product
distribution $P^n$ on $A^n$,
and any $C_n\subset A^n$,
\ben
\PR\{\mbox{``error''}\}
\;=\;1-P^n([C_n]\subD)
\; \leq \; \frac{e^{-nD^2/2}}{P^n(C_n)}.
\een
Therefore, if $\{C_n\}$ is {\em any}
$D$-covering consisting of sets
with $P^n(C_n)\approx e^{-nr}$
for some $r>0$, then the union of
the balls $B(y_1^n,D)$ centered at
the points $y_1^n\in C_n$ covers
all of $A^n$ except for a set of
probability no greater than
\be
\approx e^{-n\left(\frac{D^2}{2}-r\right)}.
\label{eq:directMC}
\ee
It is then natural to ask,
what is the {\em best} achievable error
exponent among all $D$-coverings
$\{C_n\}$ with probability no greater
that $\approx e^{-nr}$? In other words,
we are asking for small sets
with the largest possible ``boundary,''
sets $C_n$ with ``volume''
$P^n(C_n)$ no greater than $e^{-nr}$
but whose $D$-blowups $[C_n]\subD$
cover as much of $A^n$ as possible.
As pointed in \cite{kontoyiannis:sphere:01},
this question can be thought of as the
opposite of the usual isoperimetric
problem.
 
Taking $M=P$ in the general setting described above,
we obtain the answer to this question as a corollary
to our general result in the following section;
see Corollary~3.


\section{Results}

Given any $D\geq 0$ and any $R\in\RL$, 
let $E(R,D)$ denote the best achievable
error-exponent among all $D$-coverings
with mass asymptotically bounded by $2^{nR}$.
Letting ${\cal C}(R)$ denote the collection
of all  sequences of subsets $C_n$ of
$A^n$ with $\limsup_n \frac{1}{n}\log M^n(C_n)\leq R$,
define,
$$E(R,D)
\bydef
\sup_{\{C_n\}\in{\cal C}(R)}
\liminf_{n\to\infty}
	-\frac{1}{n}\log
	\Big[
	1-P^n([C_n]\subD)
	\Big],
$$
where
`log' denotes the logarithm taken to base 2.

A weaker version of this problem was recently
considered in \cite{kontoyiannis:sphere:01},
where it was shown that 
the probability of error can only decrease to
zero if $R$ is greater than $R(D;P,M)$, 
\be
R(D;P,M)
\bydef
\inf_{(X,Y):\;X\sim P,\;E\rho(X,Y)\leq D} \Big\{
        H(P_{X,Y}\|P\times P_Y)+ E[\log M(Y)] \Big\},
\label{eq:oldR}
\ee
where the infimum is taken over all jointly distributed
random variables $(X,Y)$
such that $X$ has distribution $P$
and $E\rho(X,Y)\leq D$, and $P_{X,Y}$ 
denotes the joint distribution of $X,Y$,
$P_Y$ denotes the marginal distribution
of $Y$, and $H(\mu\|\nu)$ denotes the
relative entropy between two probability
measures $\mu$ and $\nu$ on the same finite
set $S$,
$$H(\mu\|\nu)\bydef\sum_{s\in S}\mu(s)\log\frac{\mu(s)}{\nu(s)}.$$
Therefore, the error-exponent
$E(R,D)$ can only be nontrivial
(i.e., nonzero) for $R>R(D;P,M)$.
Also note that any $C_n\subset A^n$ has
$$\frac{1}{n}\log M^n(C_n)\leq \frac{1}{n}\log M^n(A^n)
=\log M(A).$$
Hence, from now on we restrict attention
to the range of interesting values for $R$
between $R(D;P,M)$ and $\Rmax\bydef\log M(A)$.

\medskip

\noindent
{\em Theorem. }
For all $D\in[0,\Dmax)$ and all $R(D;P,M)<R<\Rmax$,
the best achievable exponent of the error
probability, among all
$D$-coverings $\{C_n\}$ with mass asymptotically
bounded by $2^{nR}$, is
$$E(R,D)=E^*(R,D)
\bydef\inf_{Q\,:\,R(D;Q,M)>R} H(Q\|P),$$
where $R(D;P,M)$ is defined in (\ref{eq:oldR})
and
$H(Q\|P)$ denotes the relative entropy
(or Kullback-Leibler divergence) between
two distributions $P$ and $Q$.

\medskip

\noindent
{\em Remarks. }

{\em 1. A slightly different error-exponent. }
	Alternatively, we can define 
	a version of the optimal 
	error-exponent by considering only 
	$D$-coverings $\{C_n\}$ 
	with mass bounded 
	by $2^{nR}$ for {\em all} $n$:
	$$E'(R,D)
	\bydef
	\liminf_{n\to\infty}
        -\frac{1}{n}\log
	\left\{
	\min_{C_n\,:\,M^n(C_n)\leq 2^{nR}}
	\Big[
        1-P^n([C_n]\subD)
        \Big]
	\right\}.
	$$
	From the theorem it easily follows that
	$E'(R,D)$ is also equal to $E^*(R,D)$ at all
	points $R$ where $E^*(R,D)$ is continuous
	and, since it is nondecreasing in $R$,
	$E^*(R,D)$ is indeed continuous at all 
	except countably many values of $R$. 
	But in general it may fail to be continuous 
	everywhere, as illustrated in 
	the discussions by Marton \cite{marton:74}
	and Ahlswede~\cite{ahlswede:90} for
	the special case of lossy data compression 
	(which corresponds to taking $M(x)\equiv1$;
	see Example~2 below).


\newpage

{\em 2. Proof. }
	The proof of the theorem is a 
	modification of Marton's 
	\cite{marton:74} original argument
	for the case of error-exponents in lossy data 
	compression. The optimal 
	sets $\{C_n\}$ achieving $E^*(R,D)$ 
	are randomly generated, and they are {\em universal}
	in that their construction only depends on $R$, 
	$D$, and $M$. Therefore, they achieve the optimal 
	error-exponent simultaneously for all distributions 
	$P$.

\medskip

\noindent
{\em Example 1: Hypothesis Testing. }

Let $P_0$ and $P_1$ be two probability 
distributions $A$ with all positive probabilities. 
Suppose that the 
null hypothesis that a sample 
$X_1^n=(X_1,X_2,\ldots,X_n)$ of $n$
independent observations comes from $P_0$
is to be tested against the simple
alternative that $X_1^n$ 
comes from $P_1$. Any test between 
these two hypotheses is simply a
decision region $C_n\subset A^n$:
If $X_1^n\in C_n$ we declare that
$X_1^n\sim P_1^n$, otherwise
we declare $X_1^n\sim P_0^n$.
The set $C_n$ is called the 
{\em critical region}, and the 
type-I and type-II probabilities 
of error associated with the test
are, respectively,
\ben
\alpha_n = P_0^n(C_n)
\;\;\;\;\;
\mbox{and}
\;\;\;\;\;
\beta_n  = P_1^n(C^c_n).
\een
Clearly we wish to have
$\alpha_n$ and $\beta_n$
both decrease to zero as
fast as possible. 
In particular, we may ask 
how quickly $\beta_n$
can decay to zero if we require
that $\alpha_n$ decays 
exponentially at some 
rate $r>0$, i.e.,
$\alpha_n\approx 2^{-nr}$. 
In statistical terminology,
we are asking for the fastest
rate of decay of the type-II 
error probability
among all tests with significance
level $\alpha_n\leq 2^{-nr}$.

Formally, we want to identify the
best exponent of the
error probability
$\beta_n=1-P_1^n(C_n)$
among all $C_n$ with $P_0^n(C_n)\leq 2^{-nr}$.
Taking $P=P_1$, $M=P_0$, $R=-r$, and allowing 
{\em no} distortion, this question reduces exactly 
to the our earlier sphere-covering problem.
[To be precise, allowing no distortion means
we take $D=0$ with $\rho(x,y)$ being Hamming 
distortion as in (\ref{eq:hammingdist}).]
Accordingly, 
$R(D;P,M)=R(0;P_1,P_0)$ turns out
to be equal to $-H(P_1\|P_0)$, and from the 
theorem we immediately obtain the 
following classical result of 
Hoeffding. Also see 
\cite[Thms.~9,~10]{blahut:74} and
\cite[Ex.12,~p.43]{csiszar:book} for 
versions of this result in the 
information theory literature.

\medskip

{\em Corollary 1.} (Hypothesis Testing) \cite{hoeffding:65}
Let $\{C_n\}$ be an arbitrary sequence of tests
with associated error probabilities $\alpha_n$
and $\beta_n$ as above.
Among all tests with
$$\limsup_{n\to\infty}\frac{1}{n}\log\alpha_n\leq -r$$
for some $r\in(0,H(P_1\|P_0))$, the fastest achievable 
asymptotic rate of decay of $\beta_n$ is
$$\lim_{n\to\infty}-\frac{1}{n}\log\beta_n =
\inf_{Q\,:\,H(Q\|P_0)<r} H(Q\|P_1).$$

\medskip

As mentioned earlier, the optimal
decision regions $C_n$ in the Corollary
are randomly generated. Therefore,
although they do achieve asymptotically
optimal performance, they are not 
optimal for finite $n$
in the Neyman-Pearson sense.


\newpage

\noindent
{\em Example 2: Lossy Data Compression. }

Suppose data $X_1^n=(X_1,X_2,\ldots,X_n)$ 
is generated by a stationary, memoryless 
source, i.e., $X_1^n$ are
i.i.d.\ (independent and identically
distributed) random variables,
with distribution $P$ on the finite
alphabet $A$. The objective of lossy data
compression is to find efficient 
representations $y_1^n\in A^n$
for all source strings $x_1^n\in A^n$.
In particular,
suppose that the maximum
amount of distortion $\rho_n(x_1^n,y_1^n)$
that we are willing to tolerate between 
any source string $x_1^n$ and its 
representation $y_1^n$ is some $D\geq 0$, 
where $\{\rho_n\}$ is a family of single-letter
distortion measures as in (\ref{eq:sldistance}).
Then the problem is to find an
efficient codebook $C_n\subset A^n$
such that for most of the source strings
$x_1^n$ there is a $y_1^n\in C_n$ with
$\rho_n(x_1^n,y_1^n)\leq D$.

Here,
an efficient codebook $C_n$ is one that 
leads to good compression, i.e.,
one whose size is as small as possible. 
And, on the other hand,
we also want to make sure that the probability
that a source string cannot be represented
by any element of $C_n$ with distortion $D$ 
or less, is small.
Taking $M$ to be counting measure
($M(x)=1$ for all $x\in A$), 
the mass $M^n(C_n)$
of the codebook becomes its size $|C_n|$,
and the problem of finding a good codebook
reduces to the earlier sphere-covering question. 
Accordingly,
the rate-function $R(D;P;M)$ reduces to
Shannon's rate-distortion function $R(D;P)$,
and the theorem yields Marton's 
error-exponents result.

\medskip
 
{\em Corollary 2.} (Lossy Data Compression) \cite{marton:74}
Let $D\geq 0$ be a given distortion level, 
and $R(D;P)<R<\log |A|$.
Among all sequences of codebooks $\{C_n\}$ with asymptotic
rate no greater than $R$~bits/symbol, 
$$\limsup_{n\to\infty}\frac{1}{n}\log |C_n| \leq R,$$
the fastest achievable asymptotic rate of decay
of the probability of error is
$$\lim_{n\to\infty}-\frac{1}{n}\log
	\Big[
        1-P^n([C_n]\subD)
        \Big]
=
\inf_{Q\,:\,R(D;Q)> R} H(Q\|P).$$

%

\medskip
 
\noindent
{\em Example 3: Measure Concentration on the Binary Cube. }

Consider again the setting of the example described in the
introduction. There we asked for the {\em best} achievable 
error exponent among all $D$-coverings $\{C_n\}$ with 
probability no greater that $\approx e^{-nr}$. Taking $M=P$ 
in the theorem, we obtain the answer to this question 
in the following Corollary. 
Let $H_e(P\|Q)$
denote the relative entropy expressed
in nats rather than bits,
$H_e(P\|Q)=(\log_e 2)H(P\|Q)$,
and similarly write 
$R_e(D;P,M)=(\log_e 2)R(D;P,M)$.

\medskip
 
{\em Corollary 3.} (Converse Measure Concentration) 
Let $D\geq 0$ and $0<r<-R_e(D;P,P)$.
Among all $D$-coverings $\{C_n\}$ with 
$$\limsup_{n\to\infty}\frac{1}{n}\log_e P^n(C_n) \leq -r,$$
the fastest achievable asymptotic rate of decay
of the probability of error is
$$\lim_{n\to\infty}-\frac{1}{n}\log_e
\Big[
1-P^n([C_n]\subD)
\Big]=
\calE^*(r,D),$$
where
$$\calE^*(r,D)\bydef 
\inf_{Q\,:\,R_e(D;Q,P)> -r} H_e(Q\|P).$$
 
\medskip

Although the exponent $\calE^*(r,D)$ above is
not as explicit as $(\frac{D^2}{2}-r)$
in (\ref{eq:directMC}), it is easy to evaluate
numerically and it contains much more
useful information. For example, Figure~1 shows 
the graph of $\calE^*(r,D)$ as a function of $r$, 
for $D=0.3$, $P$ being the Bernoulli(0.4
) distribution, and $r$ running over the range 
$r\in(0.6109,0.6393)$ where $\calE^*(r,D)$
is nontrivial (i.e., finite and nonzero).
In this case, (\ref{eq:directMC}) is only useful 
when $(\frac{D^2}{2}-r)$ is positive, i.e., 
for $r\in(0,0.045)$: There (\ref{eq:directMC})
says that, whenever 
$P^n(C_n)\approx e^{-nr}$ for some
$r\in(0,0.045)$, the probability of error
decays exponentially fast. But in that range,
and in fact for all $r$ up to $\approx 0.61$,
we have $\calE^*(r,D)=\infty$ so there are sets 
$C_n$ with $P^n(C_n)\approx e^{-nr}$
and probability of error decaying 
{\em super}-exponentially fast.
Moreover, in the range $r\in(0.6109,0.6393)$ 
where $\calE^*(r,D)$ is nontrivial,
we can choose $C_n$ with 
$P^n(C_n)\approx e^{-nr}$ and 
$\PR\{\mbox{``error''}\}
\approx e^{-n\calE^*(r,D)}$.

\begin{figure}[ht]
\centerline{\epsfxsize 5.2in \epsfbox{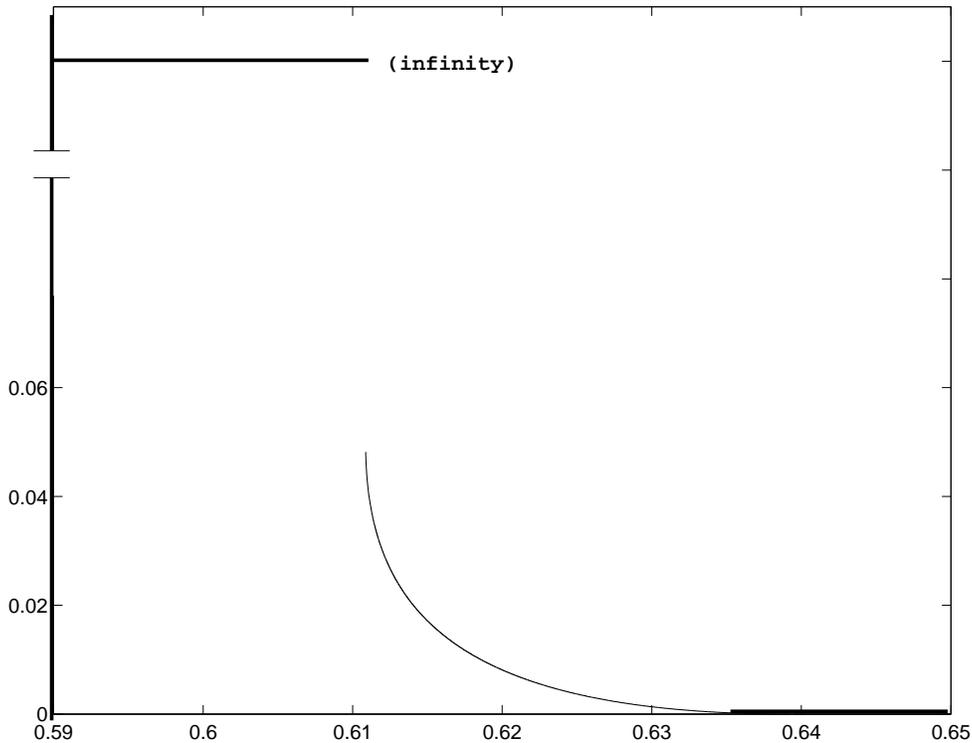}}
\caption{Graph of the error-exponent function $\calE^*(r,D)$ 
in Corollary~3 as a function of $r$, for $D=0.3$ and $P(1)=0.4$. 
Note that $\calE^*(r,D)$ is infinite for all $r\in(0,0.6109)$,
and that it is zero for $r>0.6393$.}
\end{figure}

\medskip

Finally we remark that the ``extremal'' sets 
in the classical isoperimetric problem,
namely, those $C_n$ that achieve equality in
(\ref{eq:directMC}), are very different
from the extremal sets in Corollary~3.
The former are well-known to be 
Hamming balls $B_n$ centered at
$0^n=(0,0,\ldots,0)\in A^n$,
$B_n=\{x_1^n\;:\;\rho_n(x_1^n,0^n)\leq \delta\}$
(see 
\cite{harper:66}\cite[p.~174]{mcdiarmid:89}\cite[Sec.~2.3]{talagrand:95}),
while the latter are collections
of strings $y_1^n$ randomly selected 
from a collection of suitable strings.

\medskip

 
\noindent
{\em Extensions. }
 
{\em 1. Different alphabets. }
	Although we assumed from the start
	that $\rho(x,y)$ is a distortion
	measure on $A\times A$, it is
	straightforward to generalize 
	the main result as well as the 
	subsequent discussion above to the
	case when $\rho(x,y)$ is a distortion
	measure between the ``source'' alphabet
	$A$ and a different (``reproduction'')
	alphabet $\Ahat$, as long as it is
	still the case that for each 
	$a\in A$ there exists a $b\in\Ahat$ 
	with $\rho(a,b)=0$. The necessary
	modifications to the statements and
	proofs follow exactly as in the 
	case of Marton's result; see
	\cite[Sec.~2.4]{csiszar:book}.

\medskip

{\em 2. Strong converse. } 
	As mentioned earlier, the 
	theorem is stated only for
	values of $R$ above $R(D;P,M)$
	since we trivially have $E(R,D)=0$
	for $R<R(D;P,M)$; 
	see \cite[Thm.~1]{kontoyiannis:sphere:01}.
	In that range it is also possible to prove
	a ``strong converse'' showing that,
	not only $E(R,D)=0$, but in fact
	the probability of error goes to one
	exponentially fast with a certain rate.

\section{Proof}
 
First we prove the {\em converse} 
part of the theorem, asserting that 
$E(R,D)\leq E^*(R,D)$.

Note that the rate-function $R(D;P,M)$ defined
in (\ref{eq:oldR}) is jointly uniformly continuous
in $D\geq 0$ and $P$; this can be easily seen 
to be the case by arguing along the lines of the proof of
\cite[Lemma~2.2.2]{csiszar:book}
for the rate-distortion function $R(D;P)$.
Now let $\{C_n\}$ be an arbitrary $D$-covering 
with $\{C_n\}\in{\cal C}(R)$. Take
any $Q$ on $A$ such that $R(D;Q,M) > R$ 
(if no such $Q$ exists then the claim 
is trivially true), and let $\delta>0$ 
be such that $R(D;Q,M) > R + \delta$.
Since $\{C_n\}\in{\cal C}(R)$, we have
$\log M^n(C_n)< n(R+\delta/2)$, eventually,
and by the continuity of $R(D;Q,M)$ in $D$
we can find an $\eta>0$ small enough so that
$$
\log M^n(C_n)< n(R+\delta/2)
<nR(D+\eta;Q,M),
\;\;\;\mbox{eventually.}$$

Therefore, by the ``weak converse'' in
\cite[Thm.~1]{kontoyiannis:sphere:01},
we must also have
\be
E_{Q^n}\left[\min_{y_1^n\in C_n}\rho_n(X_1^n,y_1^n)
\right]>D+\eta,
\;\;\;\mbox{eventually,}
\label{eq:converse}
\ee
where $X_1^n$ denote $n$ i.i.d.\ random
variables with distribution $Q^n$.
Writing
$$Z_n\bydef \min_{y_1^n\in C_n}\rho_n(X_1^n,y_1^n),$$
the bound in equation (\ref{eq:converse}) implies that 
$$
D+\eta<E[Z_n]\leq D\,Q^n(Z_n\leq D) + \Dmax\,Q^n(Z_n>D)$$
i.e., 
$$Q^n(Z_n>D) > \frac{\eta}{\Dmax-D}.$$

From Stein's lemma \cite[Cor.~1.1.2]{csiszar:book}
we also know that, for any $P$ and any $\epsilon>0$,
$$\lim_{n\to\infty}\frac{1}{n}
\log
\left[
\min_{ B_n \subset A^n\,:\, Q^n(B_n)> \epsilon}P^n(B_n) 
\right]= -D(Q\|P).
$$
Taking $\epsilon=\eta/(\Dmax-D)>0$ and 
applying this to the events 
$$B_n\bydef\{Z_n>D\}=[C_n]\subD^c,$$
yields
$$\liminf_{n\to\infty} 
\frac{1}{n}\log\Big[1-P^n([C_n]\subD)\Big]
\geq -D(Q\|P),
$$
and since this holds for all $Q$ with $R(D;Q,M) > R$,
we obtain
$$\limsup_{n\to\infty} 
-\frac{1}{n}\log\Big[1-P^n([C_n]\subD)\Big]
\leq E^*(R,D).$$
Finally, since $\{C_n\}\in{\cal C}(R)$ was arbitrary,
this establishes that $E(R,D)\leq E^*(R,D)$,
as required.

To prove the {\em direct} part of the theorem, asserting
the existence of a $D$-covering $\{C_n\}\in{\cal C}(R)$
such that
$$
\liminf_{n\to\infty}
        -\frac{1}{n}\log
        \Big[
        1-P^n([C_n]\subD)
        \Big]\geq E^*(R,D),
$$
we follow the same outline as in the proof of
the direct part of \cite[Thm.~2.4.5]{csiszar:book}.
 
Using the joint uniform continuity of
$R(D;P,M)$ in $D\geq 0$ and $P$,
the proof of the type-covering
lemma \cite[Lemma~2.4.1]{csiszar:book} can be
generalized to the corresponding statement with
$R(D;P,M)$ in place of $R(D;P)$. The main
new observation here is that, since all the 
elements $y_1^n$ of the covering set $B$
are drawn from the set $T^n_{[Y^*]}$ of $Y^*$-typical 
strings, where $(X^*,Y^*)$ achieve the
infimum in the definition (\ref{eq:oldR}) of
$R(D;P,M)$, their mass $M^n(y_1^n)$ satisfies
$$\frac{1}{n}\log M^n(y_1^n)\leq 
E[\log M(Y^*)] +\delta_n\left[\sum_y \log M(y)\right],$$
where the sequence $\delta_n\to 0$ as
$n\to\infty$.

Finally, following the same steps
as in the proof of the direct part of
\cite[Thm.~2.4.5]{csiszar:book} and replacing
$R(D;P)$ by $R(D;P,M)$,
we obtain the existence of a
$D$-covering $\{C_n\}\in{\cal C}(R)$
with error exponent no worse than
$E^*(R,D)-\delta,$ where $\delta>0$
is an arbitrary constant.
This proves that
$E(R,D)\geq E^*(R,D)$, and completes the
proof.
\qed

\section*{Acknowledgments}
We wish to thank Amir Dembo and Neri Merhav 
for asking us (independently) whether the 
results of \cite{kontoyiannis:sphere:01} 
could be extended to the case of 
error-exponents.


\newpage


\begin{thebibliography}{10}

\bibitem{ahlswede:90}
R.~Ahlswede.
\newblock Extremal properties of rate-distortion functions.
\newblock {\em IEEE Trans. Inform. Theory}, 36(1):166--171, 1990.

\bibitem{blahut:74}
R.E. Blahut.
\newblock Hypothesis testing and information theory.
\newblock {\em IEEE Trans. Inform. Theory}, 20(4):405--417, 1974.

\bibitem{csiszar:book}
I.~Csisz{\'{a}}r and J.~K{\"{o}}rner.
\newblock {\em Information Theory: Coding Theorems for Discrete Memoryless
  Systems}.
\newblock Academic Press, New York, 1981.

\bibitem{harper:66}
L.H. Harper.
\newblock Optimal numberings and isoperimetric problems on graphs.
\newblock {\em J. Combinatorial Theory}, 1:385--393, 1966.

\bibitem{hoeffding:65}
W.~Hoeffding.
\newblock Asymptotically optimal tests for multinomial distributions.
\newblock {\em Ann. Math. Statist.}, 36:369--408, 1965.

\bibitem{kontoyiannis:sphere:01}
I.~Kontoyiannis.
\newblock Sphere-covering, measure concentration, and source coding.
\newblock {\em IEEE Trans. Inform. Theory}, 47:1544--1552, May 2001.

\bibitem{marton:74}
K.~Marton.
\newblock Error exponent for source coding with a fidelity criterion.
\newblock {\em IEEE Trans. Inform. Theory}, 20:197--199, 1974.

\bibitem{mcdiarmid:89}
C.~McDiarmid.
\newblock On the method of bounded differences.
\newblock In {\em Surveys in combinatorics (Norwich, 1989)}, pages 148--188.
  London Math. Soc. Lecture Note Ser., 141, Cambridge Univ. Press, Cambridge,
  1989.

\bibitem{mcdiarmid:98}
C.~McDiarmid.
\newblock Concentration.
\newblock In {\em Probabilistic methods for algorithmic discrete mathematics},
  pages 195--248. Algorithms Combin., 16, Springer, Berlin, 1998.

\bibitem{talagrand:95}
M.~Talagrand.
\newblock Concentration of measure and isoperimetric inequalities in product
  spaces.
\newblock {\em Inst. Hautes \'{E}tudes Sci. Publ. Math.}, No. 81:73--205, 1995.

\end{thebibliography}

\end{document}